\documentclass[11pt,reqno]{amsart}
\usepackage{amssymb}
\usepackage[utf8]{inputenc}
\usepackage[pagewise]
{lineno}
\usepackage[margin=1.3in]{geometry}
\usepackage{amsthm}
\usepackage{mathtools}
\usepackage{tikz}
\usetikzlibrary{arrows.meta}
\usetikzlibrary{shapes.geometric}
\usetikzlibrary{patterns}
\usetikzlibrary{calc,scopes}
\usetikzlibrary{backgrounds}
\usepackage{forest}
\usepackage{wrapfig}
\usepackage{caption}
\usepackage{float}
\usepackage{units}
\usepackage{todonotes}
\usepackage{comment}
\usepackage[hidelinks]{hyperref}
\usepackage[initials,msc-links]{amsrefs}

\newcounter{indice}

\newcommand{\cyc}{\text{cyc}}

\newcommand{\cyclefig}[1]{\begin{tikzpicture}[scale= 0.35]
        \cycle{#1}
      \end{tikzpicture} \hspace{2.7mm}\vspace{0.9mm}}

\newcommand{\cycle}[1]{ 
\setcounter{indice}{0};
\foreach \i in {#1}
\addtocounter{indice}{1};
\addtocounter{indice}{1}
\draw [help lines] (0,0) grid (\theindice-1,\theindice-1);
\setcounter{indice}{1};
\foreach \i in { #1 } { 
\draw (\theindice-.5,\i-.5) [fill] circle (.18);
\draw[dashed] (\theindice-.5, \theindice-.5)--(\theindice-.5, \i-.5);
\draw[dashed] (\i-.5, \i-.5)--(\theindice-.5, \i-.5);
\addtocounter{indice}{1};
}
\addtocounter{indice}{-1};
}

\newcommand{\gridfig}[1]{\begin{tikzpicture}[scale=0.35]
    \grid{#1}
    \end{tikzpicture} \hspace{2.7mm}\vspace{0.9mm} }
    
\newcommand{\grid}[1]{
\setcounter{indice}{0};
\foreach \i in {#1}
\addtocounter{indice}{1};
\addtocounter{indice}{1}
\draw [help lines, thin] (0,0) grid (\theindice-1, \theindice-1);
\setcounter{indice}{1};
\foreach \i in {#1} {
\draw[thick] (\theindice-0.5, \theindice-0.5) -- (\theindice-0.5, \i-0.5);
\draw[thick] (\i-0.5, \i-0.5) -- (\theindice-0.5, \i-0.5);
\draw[thick] (\theindice-0.5, \i-0.5) circle (.18);
\draw[white] (\theindice-0.5, \i-0.5) [fill] circle (.12);
\draw (\i-0.5,\i-0.5) [fill] circle (.18);
\addtocounter{indice}{1};
}
\addtocounter{indice}{-1};
}

\newcommand{\xing}{\mathrm{cr}}

\newtheorem{theorem}{Theorem}[section]
\newtheorem{corollary}[theorem]{Corollary}

\newtheorem{lemma}[theorem]{Lemma}
\newtheorem{prop}[theorem]{Proposition}
\newtheorem{definition}[theorem]{Definition}

\title{Links and the Diaconis-Graham Inequality}
\author{Christopher Cornwell and Nathan McNew }

\begin{document}

\begin{abstract}
    
In 1977 Diaconis and Graham proved two inequalities relating different measures of disarray in permutations, and asked for a characterization of those permutations for which equality holds in one of these inequalities. Such a characterization was first given in 2013. Recently, another characterization was given by Woo, using a topological link in $\mathbb R^3$ that can be associated to the cycle diagram of a permutation.
We show that Woo’s characterization extends much further: for any permutation, the discrepancy in Diaconis and Graham's inequality is directly related to the Euler characteristic 
of the associated link.  This connection provides a new proof of the original result of Diaconis and Graham.  We also characterize permutations with a fixed discrepancy in terms of their associated links and find that the stabilized-interval-free permutations are precisely those whose associated links are nonsplit.

\end{abstract}

\maketitle

\section{Introduction}
\label{sec:intro}

Many different metrics can be used to measure the disarray of a shuffled list.  Treating the list as a permutation $\pi \in S_n$, and denoting by $\pi_i$ the number in position $i$ of $\pi$, some popular metrics are:
\begin{itemize}
    \item $I(\pi) = \sum_{i=1}^n \#\{j>i \mid \pi_j<\pi_i\}$, the number of inversions in $\pi$.
    \item $D(\pi) = \sum_{i=1}^n |\pi_i-i|$, the sum of the absolute  differences between the indices and values of $\pi$ (called Spearman's footrule).
    \item $T(\pi)$, the minimum number of transpositions required to produce $\pi$.
\end{itemize} 

Notice that all three metrics are zero when $\pi$ is the identity permutation. In 1849 Cayley showed that $T(\pi)$ is equal to $n$ minus the number of cycles in $\pi$, i.e.  $T(\pi) = |\pi|-\cyc(\pi)$ where $\cyc(\pi)$ is the number of cycles in $\pi$.

The so-called Diaconis-Graham inequalities \cite{diaconis} relate these three metrics for any permutation $\pi$:
\begin{equation} I(\pi) + T(\pi) \leq D(\pi) \leq 2 I(\pi). \label{eq:DGineq}
\end{equation}
Diaconis and Graham give a simple proof of the second inequality, and permutations where equality holds have a nice characterization: $D(\pi) = 2 I(\pi)$ if and only if $\pi$ contains no 3‐inversions (a triplet $i < j < k$ where $\pi_i > \pi_j > \pi_k$, equivalently an occurrence of the pattern 321). Their proof of the first inequality is not so intuitive, however, and they don't characterize when equality holds.  Such a characterization was first given in 2013 by Hadjicostas and Monico \cite{hadjicostas}.  Building on their work, and results of the present authors, Woo \cite{woo} gives a simpler, remarkable characterization of these permutations in terms of their cycle diagrams and an associated link.

The ``cycle diagram'' (or cobweb plot) of a permutation $\pi$ is obtained by plotting the permutation as follows: For each index $i$, draw a vertical line from $(i,i)$ to $(i,\pi_i)$, followed by a horizontal line to $(\pi_i,\pi_i)$. If $i$ is a fixed point of $\pi$ then there is an isolated point plotted at $(i,i)$.  Such a diagram allows a nice visualization of the cycle structure of the permutation.

\setlength{\columnsep}{15pt}%
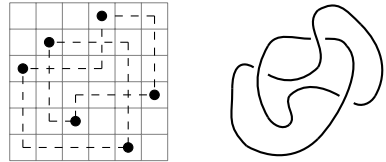
\begin{wrapfigure}{r}{0.37\textwidth}
\centering
\vspace{-3mm}
\ \ \cyclefig{4,5,2,6,1,3}\hspace{2mm}
\begin{tikzpicture}[scale=0.38]

		\node (0) at (1, 0.75) {};
		\node  (1) at (0.5, 2.75) {};
		\node  (2) at (1.75, 4.25) {};
		\node  (3) at (2.5, 1.5) {};
		\node  (4) at (3.5, 5.5) {};
		\node  (5) at (3.5, 1) {};
		\node  (6) at (5.5, 2.5) {};
		\node  (7) at (1.75, 1.75) {};
		\node  (8) at (2.5, 2.5) {};
		\node  (9) at (3.5, 3.5) {};
		\node  (10) at (4.5, 4.25) {};
		\node  (11) at (5, 5) {};
		\node  (12) at (1.25, 3.5) {};
		\node  (13) at (1.75, 3.25) {};
		\node  (14) at (3.25, 4.5) {};
		\node  (15) at (3.75, 4.5) {};
		\node  (16) at (4.25, 2.5) {};
		\node  (17) at (4.75, 2.25) {};
		\draw [thick, in=-90, out=135, looseness=1.25] (0.center) to (1.center);
		\draw [thick,in=150, out=90, looseness=1.50] (1.center) to (12.center);
		\draw [thick,in=-135, out=-30] (13.center) to (9.center);
		\draw [thick,in=-135, out=60] (9.center) to (4.center);
		\draw [thick,in=135, out=30, looseness=1.25] (4.center) to (11.center);
		\draw [thick,in=60, out=-45] (11.center) to (6.center);
		\draw [thick,in=-30, out=-120] (6.center) to (17.center);
		\draw [thick,in=45, out=150] (16.center) to (8.center);
		\draw [thick,in=45, out=-120] (8.center) to (3.center);
		\draw [thick,in=-60, out=-135] (3.center) to (7.center);
		\draw [thick,in=-120, out=120] (7.center) to (2.center);
		\draw [thick,in=-180, out=60] (2.center) to (14.center);
		\draw [thick,in=120, out=0] (15.center) to (10.center);
		\draw [thick,in=45, out=-45, looseness=0.75] (10.center) to (5.center);
		\draw [thick,in=-30, out=-135] (5.center) to (0.center);

\end{tikzpicture}
    \vspace{-2mm}
    \caption{The cycle diagram of the cycle 452613 along with the associated knot, a trefoil.}
    \label{fig:trefoil}
\end{wrapfigure}

After adding ``decorations'' at crossings, the cycle diagram may be viewed as a knot (or link) diagram. That is, considering lines in the diagram as strands of a knot, at each crossing between a vertical line and horizontal line we add the information that the vertical line is crossing over the horizontal.  Now, if each isolated point (any fixed point of $\pi$) is made to represent a small disjoint and unlinked component, then we may consider this a planar diagram of a link, which we denote $L_\pi$.\footnote{The term \emph{knot} typically requires a single connected component, but a \emph{link} may have multiple components. This happens whenever the permutation doesn't consist of a single cycle.}

The method described here of associating a link to a permutation was first introduced in \cite{cornwell-mcnew}, where it is shown that the count of the cycles of length $n$ associated to an unknotted planar diagram is enumerated by the large Schroeder numbers, and a generating function is given for the sequence enumerating permutations associated to an unlink.\footnote{An unlink is a link that is topologically the same as one in which each of the link components is a round circle, and these circles bound disks that are pairwise disjoint.}  Subsequently, Woo \cite{woo} showed that a permutation $\pi$ is associated to an unlink if and only if equality holds in the first part of \eqref{eq:DGineq}, i.e. $I(\pi)+T(\pi) = I(\pi) + n -\cyc(\pi) = D(\pi)$.  In this paper we show that this analogy extends much further.  

Denote by $\chi(L)$ the Euler characteristic of a link $L$, a well-known link invariant which is described in Section \ref{sec:Diagrams-Seifert-surfaces} below.  We also introduce the function \[x(\pi)\coloneqq D(\pi) - I(\pi) - |\pi|,\] which will be the primary permutation statistic of interest in this paper.  This quantity $x(\pi)$ is essentially the difference between the two quantities in the first inequality of \eqref{eq:DGineq}, up to the number of cycles occurring in $\pi$, i.e. $x(\pi)\ = \ D(\pi) - (I(\pi) + T(\pi)) -\cyc(\pi)$.
Our main theorem describes the relationship between these two quantities.

\begin{theorem}\label{thm:main}
    For any permutation $\pi$ the difference $x(\pi)=D(\pi) -I(\pi) - |\pi|$ is equal to the negative Euler characteristic of the associated link, i.e. \[x(\pi) = -\chi(L_\pi).\]
\end{theorem} 
  This result leads to a completely new proof of the first inequality of \eqref{eq:DGineq}, described in Section \ref{sec:corollaryDG}.

\begin{corollary}[First Diaconis-Graham inequality] \label{cor:DG}
The inequality
\begin{equation*} I(\pi) + T(\pi) \leq D(\pi)
\end{equation*}    
holds for all permutations $\pi$,
with equality holding precisely for those permutations whose diagram corresponds to an unlink (for which the Euler characteristic equals $\cyc(\pi)$).
\end{corollary}
We are able to similarly characterize those permutations $\pi$ with a given number of cycles where $D(\pi)$ exceeds $I(\pi)+T(\pi)$ by a fixed amount. The corollary below follows immediately from Theorem \ref{thm:main}.

\begin{corollary}\label{cor:exceed-by-fixed}
    For any integer $k\geq 0$, the set of permutations $\pi$ where $D(\pi)$ exceeds $I(\pi)+T(\pi)$ by exactly $k$ is precisely the set of permutations corresponding to links with Euler characteristic $\chi(L_\pi)=-k+\cyc(\pi)$.
\end{corollary}%

Before giving the proof of Theorem \ref{thm:main} in Section \ref{sec:mainproof}, we prove several lemmata in Section \ref{sec:lemmata} that will be needed in the proof.  From these results we obtain the following theorem, which may be of independent interest.  

A link $L$ in $\mathbb R^3$ is called a split link if there exists a surface in $\mathbb R^3$, homeomorphic to a 2-sphere, which is disjoint from $L$ and has components of $L$ in both its interior and exterior.  If no such surface exists, $L$ is called non-split.  We are able to characterize those permutations corresponding to non-split links as follows.

\begin{theorem} \label{thm:nonsplit}
    A permutation $\pi$ corresponds to a non-split link $L_\pi$ if and only if $\pi$ is a stabilized-interval-free permutation.
\end{theorem}

The stabilized-interval-free permutations were defined and enumeratated by Callan \cite{callan}. See Section \ref{subsec:nonsplit} for details.  The following corollary then follows from his enumeration.

\begin{corollary} \label{cor:nonsplit}
    The probability a randomly chosen derangement of $n$ corresponds to a non-split link tends to 1 as $n\to \infty$.
\end{corollary}

\section{Cycle Diagrams and Seifert Surfaces of Links}
\label{sec:Diagrams-Seifert-surfaces}

As discussed in the introduction, the cycle diagram of a permutation $\pi$ may be decorated with crossing information and viewed as a link diagram. We write $G_\pi$ for this decorated cycle diagram and $L_\pi$ for the associated link in $\mathbb R^3$. If there are no fixed points of $\pi$, the diagram $G_\pi$ is referred to in the literature as a grid diagram of the link (see \cite{ng-thurston} for a survey). 

\begin{wrapfigure}{r}{0.33\textwidth}
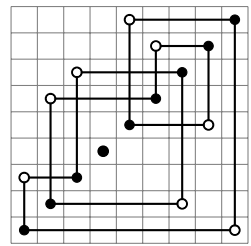

\centering
\vspace{-4.5mm}
\hspace{1mm}\gridfig{3,6,7,4,9,8,2,5,1}
\caption{The diagram $G_\pi$, where $\pi=367498251$.}
    \label{fig:467513298}
\end{wrapfigure}

Standard definitions of grid diagrams do not allow for isolated points, such as $G_\pi$ might contain on the diagonal. However, a grid diagram for $L_\pi$ is easily obtained by replacing each isolated point by a small square contained in adjacent rows and columns, thereby producing a grid diagram with size $|\pi|+\#\{\text{fixed points}\}$. 

We orient the strands of $L_\pi$ so that vertical lines are oriented away from the diagonal; tracing the diagram in the order that $\pi$ permutes its elements agrees with this orientation of the link.

Each crossing in an oriented link diagram has a sign, positive or negative. Our orientation convention causes every crossing in $G_\pi$ to be negative, so every link $L_\pi$ is the mirror of a \emph{positive link}. In addition, if $L_1$ and $L_2$ are two link components of $L_{\pi}$ then the \emph{linking number} (cf.~\cite[Chapter 5]{Rolfsen}) of $L_1$ and $L_2$ must be non-positive, and is only 0 if the planar projections of $L_1$ and $L_2$ into $G_\pi$ are disjoint. 

\begin{prop}\label{prop:split-links}
    For a permutation $\pi$, the link $L_{\pi}$ is a split link if and only if $G_{\pi}$ is disconnected (as a subset of the plane).
\end{prop}
\begin{proof}
The linking numbers of pairs of link components are a link invariant. Thus, given a \emph{split link}, for any of its components $L_1$ that is on one side of the splitting sphere, $L_1$ has linking number 0 with every component found on the other side of the sphere. Hence, using the previous observation, there is a proper subset of the set of link components of $L_\pi$ whose planar projection in $G_{\pi}$ is disjoint from the projection of link components that are not in this subset (i.e., $G_\pi$ is disconnected as a union of lines and points in the plane).
\end{proof}

Write $\xing(G_\pi)$ for the number of crossings in $G_\pi$. Another statistic of $G_\pi$ of interest is the number of upper-right corners in the diagram: indices $i$, $1\le i\le n$, which satisfy both $\pi(i)<i$ and $\pi^{-1}(i) < i$. Define $s(G_\pi)$ to be the number of upper-right corners plus the number of fixed points of $\pi$. (Note: each small square replacing an isolated diagonal point would have exactly one upper-right corner.) 

Any link in $\mathbb R^3$ may be regarded as the boundary of an oriented surface  embedded in $\mathbb R^3$, a \emph{Seifert surface} of the link. Given a diagram of the link (with orientation), Seifert described an algorithm to determine such a surface (\cite{Seifert}, cf.\ \cite[Chapter 5]{Rolfsen}). Our manner of proving the main result is informed by this algorithm, so we describe it here.

Given an oriented link diagram $Z$, remove each crossing in $Z$ by changing how the strands are connected in a neighborhood of the crossing, as in Figure \ref{fig:remove-crossing} below. Note that the change is made to be coherent with the orientation on strands outside of the neighborhood. The result $\hat Z$ is a union of oriented simple closed curves, no two of which intersect, called \emph{Seifert circles}. As we will use them again in the algorithm, we remember the arcs on Seifert circles that were paired inside one of the neighborhoods.

\begin{figure}[H]
    \centering
	\begin{tikzpicture}[>=stealth, scale=1.2]
	\draw[->]   (-1,0)--(0,1);
	\draw[white,line width=3]     (0,0)--(-1,1);
	\draw[->]   (0,0)--(-1,1);
    \draw[gray] (-0.5,0.5) circle (20pt);
 
	\draw[dashed,line width=2, ->] (0.5,0.5) -- (1,0.5);
	
	\draw[->] 	(1.5,0)..controls (1.75,0.4) and (1.75,0.6)..(1.5,1);
	\draw[->]	(2.5,0)..controls (2.25,0.4) and (2.25,0.6)..(2.5,1);
    \draw[gray] (2,0.5) circle (20pt);
	\end{tikzpicture}
\caption{Removing a crossing}
\label{fig:remove-crossing}
\end{figure}
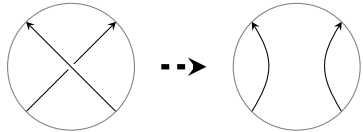

Each Seifert circle bounds a disk, oriented so that the orientation on the Seifert circle agrees with the induced boundary orientation. Place the disks in $\mathbb R^3$ at heights so that, if Seifert circles $S$ and $S'$ are such that $S$ is in the interior of $S'$ (in the plane of the diagram), then the disk for $S$ is above the disk for $S'$. 

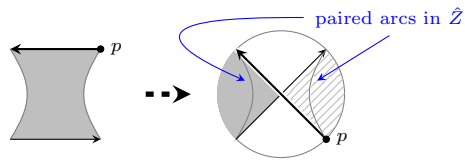
\begin{figure}[H]
    \centering
	\ 
 \hspace{1cm}\begin{tikzpicture}[>=stealth, scale=1.2]
    \filldraw[gray!50!white] (-1,0) -- (0,0) ..controls (-0.25,0.4) and (-0.25,0.6)..(0,1)
         -- (-1,1)
        ..controls (-0.75,0.6) and (-0.75,0.4)..(-1,0);
    \draw[gray] (0,0) ..controls (-0.25,0.4) and (-0.25,0.6)..(0,1);
    \draw[gray] (-1,1) ..controls (-0.75,0.6) and (-0.75,0.4)..(-1,0);
    \draw[->] (-1,0) -- (0,0);
    \draw[->, thick] (0,1) -- (-1,1);
    \filldraw (0,1) circle (1pt);
    \draw (0,1) node[right]{{\tiny $p$}};
 
	\draw[dashed,line width=2, ->] (0.5,0.5) -- (1,0.5);
	
	\filldraw[gray!50!white]   (1.5,1) -- (2,0.5) -- (1.5,0)	
        ..controls (1.75,0.4) and (1.75,0.6)..(1.5,1);
	\filldraw[draw=none,pattern=north east lines,pattern color=gray!50!white]   (2.5,1) -- (2,0.5) -- (2.5,0)
        ..controls (2.25,0.4) and (2.25,0.6)..(2.5,1);
    \draw[gray] (2,0.5) circle (20pt);
    \filldraw[draw=none,pattern=north east lines,pattern color=gray!50!white] ([shift=(-45:20pt)]2,0.5) arc (-45:45:20pt) 
            (2.5,1)..controls (2.27,0.6) and (2.27,0.4)..(2.5,0);
	\filldraw[gray!50!white] ([shift=(135:20pt)]2,0.5) arc (135:225:20pt) 
            (1.5,0)..controls (1.75,0.4) and (1.75,0.6)..(1.5,1);
	\draw[gray] 	(1.5,0)..controls (1.75,0.4) and (1.75,0.6)..(1.5,1);
	\draw[->, thin]   (1.5,0)--(2.5,1);
	\draw[white,line width=2]     (2.5,0)--(1.5,1);
	\draw[->, thick]   (2.5,0)-- (1.5,1);
    \draw[gray]	(2.5,0)..controls (2.25,0.4) and (2.25,0.6)..(2.5,1);
    \filldraw (2.5,0) circle (1pt);
    \draw (2.5,0) node[right]{{\tiny $p$}};
    \draw[blue] (3.2,1.35) node {{\tiny paired arcs in $\hat Z$}};
    \draw[->,blue] (2.25,1.35) ..controls (1.5,1.35) and (0,1.35) ..(1.6,0.65);
    \draw[->,blue] (3.2,1.15) -- (2.4,0.65);
    \end{tikzpicture}
\caption{A band is attached with a half-twist; diagonal pattern represents the opposite side of the surface.}
\label{fig:twisted-band}
\end{figure}

Finally, for each crossing from $Z$ we consider the associated arcs in $\hat Z$ that were paired. Viewing these arcs as part of the boundary of disks for Seifert circles, attach at those arcs two opposite sides of a half-twisted band, as in Figure \ref{fig:twisted-band} (with the remaining, unattached, part of the band's boundary now agreeing with the original crossing strands from $Z$). This completes the construction of the surface.

In Seifert's original algorithm, if the surface is disconnected after attachment of half-twisted bands then handles are added to connect the surface. This is needed when the link diagram is disconnected as a subset of the plane. In this paper, we leave the surface disconnected in this case.\footnote{For the sake of clarity, we highlight the difference between components of a diagram $G_{\pi}$ and the link components of $L_{\pi} \subset \mathbb R^3$. For example, the diagrams depicted in Figure \ref{fig:example_translation} each have a single connected component; however, the associated link has 3 link components.}

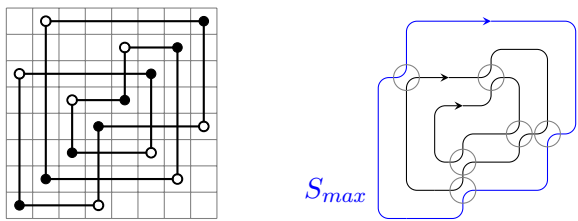
\begin{figure}[ht]
    \centering
	\gridfig{6,8,5,1,7,3,2,4} 
	\hspace{0.5cm} \begin{tikzpicture}[>=stealth,scale=3.75]
	\draw [rounded corners, blue, ->] (0.1,0)--(0,0)--(0,0.5)--(0.1,0.5)--(0.1,0.7) -- (0.4,0.7);
    \draw [rounded corners, blue] (0.4,0.7)--(0.7,0.7)--(0.7,0.3)--(0.6,0.3)--(0.6,0.1)--(0.3,0.1)--(0.3,0)--(0.1,0);
	\draw [rounded corners, ->] (0.2,0.1) -- (0.1,0.1)--(0.1,0.5) -- (0.25,0.5);
    \draw [rounded corners] (0.25,0.5)--(0.4,0.5)--(0.4,0.6)--(0.6,0.6)--(0.6,0.3)--(0.5,0.3)--(0.5,0.2)--(0.3,0.2)--(0.3,0.1)--(0.2,0.1);
	\draw [rounded corners, ->] (0.25,0.2)--(0.2,0.2)--(0.2,0.4)--(0.3,0.4);
    \draw [rounded corners] (0.3,0.4)--(0.4,0.4)--(0.4,0.5)--(0.5,0.5)--(0.5,0.3)--(0.3,0.3)--(0.3,0.2)--(0.25,0.2);
	\draw [blue] (-0.15, 0.1) node {$S_{max}$};

    \draw[gray] (0.1,0.5) circle (1.3pt);
    \draw[gray] (0.3,0.2) circle (1.3pt);
    \draw[gray] (0.3,0.1) circle (1.3pt);
    \draw[gray] (0.4,0.5) circle (1.3pt);
    \draw[gray] (0.5,0.3) circle (1.3pt);
    \draw[gray] (0.6,0.3) circle (1.3pt);
	\end{tikzpicture}
\caption{The grid diagram obtained from the permutation 68517324 and the Seifert circles obtained from the link diagram, with orientation indicated by arrows.  The circles surround each of the removed crossings.}
\label{fig:Seifert-circles}
\end{figure}

Note that for a permutation $\pi$, the orientation convention and the fact that vertical segments cross over horizontal in $G_\pi$ make the result of removing a crossing always appear the same (see Figure \ref{fig:Seifert-circles}). A consequence is that, for any Seifert circle of $G_\pi$, above the diagonal it is (weakly) increasing in the vertical and horizontal directions; below the diagonal it is decreasing in the vertical and horizontal directions. Thus, there is exactly one upper-right corner and one lower-left corner on each Seifert circle. That is, $s(G_\pi)$ equals the number of Seifert circles from the algorithm.

Given a permutation $\pi$, we let $F_\pi$ denote the surface constructed from the grid diagram $G_\pi$ through the algorithm above. The main result presented in this paper relates the value of the permutation statistic $x(\pi) = D(\pi)-I(\pi)-|\pi|$ to a numerical invariant of the surface $F_\pi$, namely the \emph{Euler characteristic} of $F_\pi$. 

The Euler characteristic of a surface $F$ may be defined via a triangulation of the surface (see \cite[pg. 85{--}87]{Shifrin} for a brief, but informative, introduction to triangulations of surfaces). If a triangulation of $F$ has $v$ vertices, $e$ edges, and $f$ faces, the Euler characteristic is $\chi(F) = v - e + f$. While the Euler characteristic of $F$ may be computed this way, it is an invariant of $F$ and is independent of choice of triangulation.

To define the Euler characteristic of a link, we consider orientable surfaces having the link as boundary. Were we to allow for any such surface, then the Euler characteristic could be arbitrarily increased by adding a separate surface component that is a sphere (and has no boundary). To avoid this pitfall, but still allow for disconnected surfaces, we require that every surface component have a non-empty boundary.
\begin{definition}
    We say that an orientable surface $F$ is a (potentially disconnected) \emph{Seifert surface} of an oriented link $L$ if every component of $F$ has non-empty boundary and the boundary of  $F$ (with induced orientation) agrees with $L$. The \emph{Euler characteristic of }$L$, written $\chi(L)$ is the maximal Euler characteristic of a Seifert surface of $L$.
\end{definition}

A triangulation of $F_\pi$ may be chosen so that there is an edge corresponding to each paired arc where the half-twisted bands were attached to a disk of a Seifert circle. 
Note that a disk has Euler characteristic equal to 1. By choosing the triangulation to be such that each of the half-twisted bands consists of exactly two triangles in the triangulation (sharing one ``diagonal'' edge in the interior of the band), one may check that  
\begin{equation} \label{eq:chi} \chi(F_\pi) = s(G_\pi) - \xing(G_\pi).\end{equation}

There is a relationship between grid diagrams and Legendrian links that is often exploited in the knot theory (and Legendrian knot theory) literature. We use this relationship to see that $F_\pi$ has maximum Euler characteristic among Seifert surfaces of $L_\pi$. 

\begin{prop}\label{prop:Euler-characteristic}
If $F$ is a Seifert surface of $L_\pi$, then $\chi(F) \le \chi(F_{\pi})$, and so $\chi(L_\pi) = \chi(F_\pi)$.
\end{prop}
   
\begin{proof}Let $G_{\pi}^{\circlearrowright}$ be the grid diagram obtained by rotating $G_\pi$ clockwise 90$^\circ$, changing all crossings to retain that vertical strands cross over horizontal. There is a \emph{Legendrian} link $\Lambda_\pi$ associated to $G_\pi^{\circlearrowright}$ which has the topological type of the mirror of $L_\pi$ (see \cite{ng-thurston}). One of the classical invariants of a Legendrian link $\Lambda$ is the Thurston-Bennequin number $\mathrm{tb}(\Lambda)$. For the Legendrian $\Lambda_\pi$, its value comes out as $\mathrm{tb}(\Lambda_\pi) = \xing(G_\pi) - s(G_\pi)$. Now, by the Eliashberg-Bennequin inequality~\cite{eliashberg}, we have $\mathrm{tb}(\Lambda_\pi) \le -\chi(F)$ for any surface $F$ with $\Lambda_\pi$ as its boundary link (not allowing connected components of $F$ to have empty boundary). Since mirroring leaves the Euler characteristic of a surface unchanged, and noting that $-\chi(F_\pi) = \xing(G_\pi) - s(G_\pi)$, we obtain the desired result.
\end{proof}

\section{Proof of Corollary \ref{cor:DG}}
\label{sec:corollaryDG}

From the classification of surfaces, a connected surface with one boundary component has Euler characteristic at most $1$ (and a connected surface with more than one boundary component has Euler characteristic less than $1$). Let $F$ be an oriented surface with $b$ boundary components, such that every surface component has non-empty boundary. Since $\chi(F)$ equals the sum of the Euler characteristics of the connected components of $F$, we have that $\chi(F) \le b$. 
Given a permutation $\pi$, if $b$ is the number of cycles of $\pi$ then $b$ is also the number of link components in $L_\pi$. This provides a new proof of the Diaconis-Graham inequality.
\begin{corollary}Given any permutation $\pi$, $$   I(\pi)+T(\pi) \ = \ I(\pi) + |\pi|-\mathrm{cyc}(\pi) \ \le  \ D(\pi).$$
\end{corollary}
\begin{proof}By Theorem \ref{thm:main}, \[D(\pi) - I(\pi)-|\pi|+\mathrm{cyc}(\pi) = -\chi(F_\pi)+b \ge 0 \]
since $b = \cyc(\pi)$ counts both the number of boundary components of the surface $F_\pi$ and the number of cycles appearing in $\pi$. 
\end{proof}

We also obtain a different proof of Alex Woo's \cite{woo} characterization of the permutations associated to an unlink.
\begin{corollary}
    A permutation is associated to an unlink if and only if it is \emph{shallow} as defined in \cite{berman}, namely $D(\pi) - I(\pi) -|\pi| +\cyc(\pi) = 0$.
\end{corollary}  This follows directly from Theorem \ref{thm:main}. In particular, our Theorem implies that a permutation is shallow if and only if $-\chi(L_\pi) = -b$; this occurs exactly when each component of $F_\pi$ is a disk, implying $L_\pi$ is an unlink.  Unlike the proof in \cite{woo} this approach doesn't rely on the original recursive characterization of shallow permutations given by Hadjicostas and Monico. 

\section{Lemmata} \label{sec:lemmata}
We define the direct sum of two permutations, $\sigma \oplus \tau$ in the usual way, if $\sigma$ is a permutation of length $m_1$ and $\tau$ is a permutation of length $m_2$ then $\sigma \oplus \tau $ is a permutation of $m_1+m_2$ given by
\[(\sigma \oplus \tau)(i) = \begin{cases}
\sigma(i) & 1\leq i \leq m_1 \\
\tau(i-m_1) + m_1 & m_1 < i \leq m_1+m_2. \end{cases}\]
\begin{lemma} \label{lem:directsum} The direct sum of two permutations satisfies both \[x(\sigma \oplus \tau) = x(\sigma)+x(\tau)\]
    and   \[ \chi(L_{\sigma \oplus \tau}) = \chi(L_\sigma)+\chi(L_\tau).\]
\end{lemma}

\begin{proof}
    It follows immediately from the definitions that both $D(\sigma \oplus \tau) = D(\sigma)+D(\tau)$ and $I(\sigma \oplus \tau) = I(\sigma)+I(\tau)$, and so the same holds for the function $x$.  Since the component of the diagram corresponding to $\sigma$ is disconnected from the component corresponding to $\tau$, and our convention is to construct disconnected surfaces for each component of the diagram of $G_{\sigma \oplus \tau}$, we find that 
    \[ \chi(L_{\sigma \oplus \tau}) = \chi(F_{\sigma \oplus \tau}) = \chi(F_{\sigma}) + \chi(F_\tau)= \chi(L_\sigma)+\chi(L_\tau). \qedhere \]
\end{proof}

The following lemma shows that a ``translation'' of a grid diagram and corresponding permutation leaves both the link type and the statistic $x(\pi)$ unchanged.  

\begin{lemma} \label{lem:translation}
    For a permutation $\pi \in S_n$, and the permutation (in one line notation) $\sigma = 234\cdots n1 \in S_n$ (note $\sigma = (123 \cdots n)$ in cycle notation) define $\pi' = \sigma \pi \sigma^{-1}$. Then the link type of $L_\pi$ is the same as that of $L_{\pi'}$ and also $$x(\pi)=x(\pi').$$
\end{lemma}

\begin{figure}[ht]
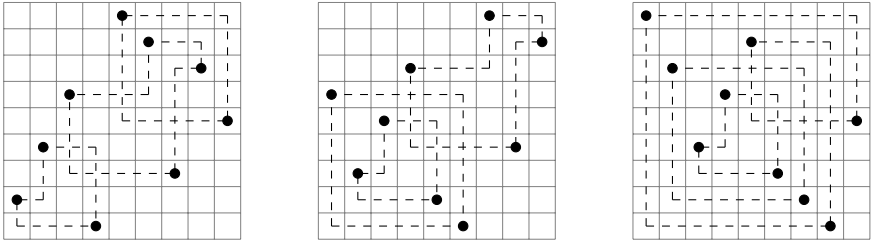

	\centering
	\ \cyclefig{2,4,6,1,9,8,3,7,5}
	\hspace{0.5cm} 
    \cyclefig{6,3,5,7,2,1,9,4,8}
    \hspace{0.5cm} 
    \cyclefig{9,7,4,6,8,3,2,1,5}
	
	\caption{Successive applications of a horizontal and vertical translation of the permutation and diagram obtained from $\pi=246198375$ (resulting in $\pi'=635721948$ and $\pi''=974683215$).  While the diagrams look relatively different, one can check that each corresponds to the same link,  and we find that $x(\pi) \ = \ D(\pi)-I(\pi)-|\pi| \ =\ 24  - 14 - 9 \ = \ 26 - 16 -9 = 38 - 28 - 9 = 1$ respectively in each of the figures above.\label{fig:example_translation}}
\end{figure}

In the terminology of grid diagrams, the grid diagrams $G_\pi$ and $G_{\pi'}$ are related by two \emph{translations} (one being horizontal and the other being vertical).  By viewing the grid diagram as being on a torus, identifying opposite edges of the square, the conjugation action has the result of translating each column to the right by one, and each row up by one. The original rightmost column (respectively topmost row) appears on the far left side (bottom side) of the grid afterwards (see \cite{ng-thurston} for details and Figure \ref{fig:example_translation} for an example).  

\begin{proof}[Proof of Lemma \ref{lem:translation}]
It is well known (\cite{ng-thurston}) that the translation operation on a grid diagram leaves the associated link type unchanged. It then follows that $L_\pi$ and $L_{\pi'}$ have the same link type, and 
\begin{equation}
-\chi(L_\pi) \ = \ -\chi(L_{\pi'}).
\end{equation}
Thus, we need only to show that $x(\pi)=x(\pi')$.  Write out $\pi = \pi_1\pi_2\ldots\pi_n$ in one-line notation and define $a, b$ to be the integers such that $\pi_a = n$ and $\pi_n = b$. We can express the permutation $\pi'$ as $\pi' = \pi'_1\pi'_2\ldots\pi'_n$, where
$$\begin{cases} 
    \pi'_{k+1} = \pi_k + 1,\qquad\text{for } 1\le k\le n-1, k\ne a \\
    \pi'_{a+1} = 1 \\
    \pi'_{1} = b+1.
\end{cases}$$
In the edge case that $a=n=b$, the first line, indicating the value of $\pi_{k+1}'$, $1\le k\le n-1$, is unchanged; however, instead of the second and third line we have that $\pi_1' = 1$.

Now we may compute that 
\begin{align*}
    D(\pi') = \sum_{i=1}^n|\pi'_i - i|    &= a + b + \sum_{\substack{1\le i<n\\ i\ne a}} |\pi_i+1 - (i+1)| \\
    &= D(\pi) + 2a + 2b - 2n
\end{align*}
in the case when $a\ne n$. If it is the case that $a=n$ (and so $b=n$ also) then \[D(\pi') = \sum_{i=2}^n|\pi'_i - i| = \sum_{i=1}^{n-1}|\pi_i - i| = D(\pi),\] since $|\pi'_1 - 1|$ and $|\pi_n-n|$ are both zero.

In order to understand how the inversions change between $\pi$ and $\pi'$, note that if $i,j \not\in \{a,n\}$, then $(\pi_i, \pi_j)$ are in the same relative position to each other in $\pi$ (i.e.\ whether $\pi_i$ is \emph{left of} or \emph{right of} $\pi_j$ in one-line notation), as the relative position of $(\pi'_{i+1}, \pi'_{j+1})$ in $\pi'$.

Since every $\pi_i \ne b$ appears to the left of $\pi_n = b$ in $\pi$, we get that the two elements $\pi_i, \pi_n$ make an inversion if and only if $\pi_i > b$, a total of $(n-b)$ inversions in $\pi$ that involve $\pi_n$. After translation, $\pi'_1=b+1$ is to the left of every other element in $\pi'$, so it makes an inversion with $\pi'_{i+1}$ if and only if $b + 1 > \pi'_{i+1}$, a total of $b$ inversions in $\pi'$ involving $\pi'_1$ (or $0$ inversions in the case that $b=n$, since $\pi_1'=1$ in that case). Thus, if $b\ne n$, then in going from $\pi$ to $\pi'$ there is a net change of $b - (n-b)=2b-n$ inversions involving the element $\pi_n$ (in $\pi$) or $\pi'_1$ (in $\pi'$). In the edge case $b=n$, we get that there is a net change of $0$ inversions involving $n$ (in $\pi$) and $0$ involving $1$ in $\pi'$.

We can make similar considerations for inversions that involve $\pi_a$ in $\pi$ and $\pi'_{a+1}=1$ in $\pi'$ (assuming $a\ne n$). One of these, already considered above, is the inversion involving $\pi_a, \pi_n$ that associate with the inversion in $\pi'$ involving $\pi'_1, \pi'_{a+1}$. The net change in inversions going from $\pi$ to $\pi'$ involving the elements $\pi_a$ or $\pi'_{a+1}$ respectively is $a - (n-a) = 2a - n$. 

When $a\ne n$ and $b\ne n$, we have shown that
\begin{align*}
    D(\pi') - I(\pi')    &= D(\pi) + 2a + 2b - 2n - (I(\pi) + 2b-n + 2a - n) \\ 
                                              &= D(\pi) - I(\pi).
\end{align*}
In the case that $a = n$ (forcing $b=n$ also), we have both $D(\pi) = D(\pi')$ and $I(\pi) = I(\pi')$.  Hence, we obtain the desired equality $x(\pi') = x(\pi)$.
\end{proof}

\begin{lemma} \label{lem:Seifertcircles}
    The Seifert circles obtained from a connected component of a cycle diagram admit a partial order under containment, and there is a unique maximal Seifert circle $S_{\text{max}}$ containing all of the others.  Furthermore, each Seifert circle has a unique upper-right and lower-left corner, and both are on the main-diagonal.
\end{lemma}

\begin{proof}
    These results are all demonstrated in the proofs of Lemma 20 and Lemma 22 (and remark 21) of \cite{cornwell-mcnew}.  While the statements of these lemmas refer only to diagrams corresponding to cycles, it is straightforward to see that all of the arguments go through identically for any connected component of a diagram.
\end{proof}

We say that a permutation $\pi$ \textit{stabilizes} an interval $I$ if $\pi(i) \in I$ for every $i\in I$. 
\begin{lemma} \label{lem:stab}
    The diagram $G_\pi$ of a permutation $\pi$ is disconnected (i.e. the union of any fixed points, vertical and horizontal segments of the diagram is disconnected as a subset of the plane) if and only if the permutation $\pi$ stabilizes a nontrivial subinterval $I\subset [1,n]$.
\end{lemma}
\begin{proof}
It follows immediately that if $\pi$ stabilizes an interval $I$, then a square drawn to encompass the points on the diagonal corresponding to the integers in $I$ will disconnect the components of the diagram of $\pi$. So we focus on the reverse direction.  

Suppose that the diagram of $\pi$ is disconnected. This occurs if and only if there is a simple closed curve in the plane, not intersecting any point on $G_\pi$, which separates the diagram into two (non-empty) \emph{sub-diagrams} {--} one in the interior, the other in the exterior of the closed curve. Call these sub-diagrams $G^{(1)}_\pi$ and $G^{(2)}_\pi$.  Note that the diagonal nodes $\{1,\ldots,|\pi|\}$ are split into two subsets $s^{(1)}(\pi)$ and $s^{(2)}(\pi)$ (where $i\in s^{(k)}(\pi)$ if and only if diagonal node $i$ is part of $G^{(k)}_{\pi}$). 

Assume that $1 \in s^{(1)}_\pi$.  If $s^{(1)}_\pi$ corresponds to an interval we are done, so suppose it does not, and let $n_1,n_2 \in s^{(1)}_\pi$ be the smallest pair of integers such that $n_2 > n_1+1$ and the entire interval $I=[n_1+1,n_2-1]\subseteq s_\pi^{(2)}$.  If $\pi$ stabilizes $I$, then again we are done.  So we now suppose it does not and derive a contradiction, which will complete the proof.  

If $\pi$ does not stabilize $I$, then for some $i\in I$, $\pi(i) \notin I$, and hence, $ \pi(i)>n_2$ (by the minimality of $n_1,n_2$).  Note that since the vertical line segment in $G^{(2)}_\pi$ that extends upward from $i$, extends above height $n_2$, we have $\pi^{-1}(n_2)>i$ (otherwise, $i$ and $n_2$ could not be in disconnected sub-diagrams, as the horizontal segment of $G^{(1)}_{\pi}$ that is at height $n_2$ would intersect the vertical one extending upward from $i$). Thus $\pi^{-1}(n_2)>n_2$.  Similarly, there must be $j>n_2$ with $\pi(j)\in I$, and so, in order to avoid a crossing it must be that $\pi(n_2)>n_2$.  Thus, we find that the line segments adjacent to the point $(n_2,n_2)$ form a lower-left corner of some Seifert circle $S_1$.  

Appealing to Lemma \ref{lem:Seifertcircles}, each Seifert circle has a unique lower-left and upper-right corner (both of which lie on the main diagonal). We find that $S_1$ must lie completely inside the Seifert circle $S_2$ that contains the line segments connected to the point $(i,\pi_i)$ (whose lower left corner must have $x$-coordinate at most $i<n_2$, and whose upper right corner has $y$-coordinate at least $\pi_i>n_2$).  Since the line segments connected to $(i,\pi_i)$ are part of $G_\pi^{(2)}$, all of the line segments that were used to create $S_2$ must also be part of $G_\pi^{(2)}$. However, they completely surround the point $(n_2,n_2)$, disconnecting it from the diagonal points at or below $n_1$, which is not possible with $G_\pi^{(1)}$ and $G_\pi^{(2)}$ being separated by a simple closed curve.
\end{proof}

\subsection{Non-split links and stabilized-interval-free permutations}  \label{subsec:nonsplit}

Callan \cite{callan} defined a permutation $\pi$ to be  \textit{stabilized-interval-free} if there is no nontrivial interval $I \subset \{1,2,\ldots |\pi|\}$ stabilized by $\pi$.  He then showed that such permutations are enumerated by sequence \href{https://oeis.org/A075834}{A075834} in the OEIS \cite{oeis}, the integer sequence whose generating function $A(x)$ has the property that the coefficient of $x^{n-1}$ in $A(x)^n$ is $n!$.  

Theorem \ref{thm:nonsplit}, that the permutations whose diagrams correspond to non-split links are precisely the stabilized-interval-free permutations, follows from Lemma \ref{lem:stab} and Proposition \ref{prop:split-links}.

Since the $n$-th term in the sequence \href{https://oeis.org/A075834}{A075834} is asymptotically $\frac{n!}{e} \left(1 - 1/n +O(n^{-2})\right)$, which is almost the same as that of the derangements (whose count is $\frac{n!}{e} +O(1)$), and of which the stabilized-interval-free permutations are a subset, we find that the probability a randomly chosen derangement of $n$ corresponds to a non-split link tends to 1 as $n\to \infty$, as noted in Corollary \ref{cor:nonsplit}.

\section{Proof of Theorem \ref{thm:main}} \label{sec:mainproof}
We will prove that $x(\pi) = -\chi(L_\pi)$ by induction on the length of the permutation, $|\pi|$.  An interesting problem would be to find a non-inductive proof of this theorem.
 
 For permutations of length at most 2, the result is immediate: in every case the corresponding link $L_\pi$ is an unlink. (When $\pi=1$ or $\pi=21$, the permutation has a single cycle, the link has a single component, and $x(\pi)=-\chi(L_\pi)=-1$. When $\pi=12$ the permutation and link have 2 cycles/components, and $x(\pi)=-\chi(L_\pi)=-2$.

We now suppose that the result holds for permutations of all lengths less than $n$ and that $|\pi| = n>2$.  For any permutation of length $n$ one of the following cases hold:

\begin{itemize}
    \item Case 1: The cycle diagram $G_\pi$ is disconnected.
    \item Case 2: The diagram $G_\pi$ is connected but there exists an index $1<i<n$ such that either $\pi_j\leq i$ for all $j<i$ or $\pi_j\geq i$ for all $j>i$. (In this situation the outermost Seifert circle touches the main diagonal at the index $i$.)
    \item Case 3: The diagram $G_\pi$ is connected and no index with the property described in case 2 exists. (The outermost Seifert circle doesn't touch the main diagonal, except at the bottom left and upper right corners.)
\end{itemize}
In each case we will see that the result can be reduced to the case of permutations of smaller length, which will complete the induction.

\subsection{Case 1.} \textit{$G_\pi$ is disconnected.} 
\smallskip

In this case, by Lemma \ref{lem:stab}, there exists a nontrivial subinterval $I\subset [1,n]$ so that $\pi$ stabilizes $I$.  By repeated application of Lemma \ref{lem:translation} (if necessary) we can translate the grid diagram of $\pi$ to the grid diagram of another permutation $\pi'$ (with $x(\pi)=x(\pi')$ and $-\chi(L_\pi) = -\chi(L_{\pi'})$) such that the interval $I$ has been translated to an interval, stabilized by $\pi'$, of the form $[a,n]$ for some integer $a$.  Now, it is clear that $\pi'$ also stabilizes the complementary interval $[1,a-1]$ and so $\pi'$ decomposes as a direct sum $\pi'=\rho\oplus \tau$.  Since both $\rho$ and $\tau$ are permutations of length shorter than $|\pi|$, we are done by Lemma \ref{lem:directsum} and induction.

\subsection{Case 2.} \textit{The diagram $G_\pi$ is connected and there is an index $1<i<n$ such that either $\pi_j\leq i$ for all $j<i$ or $\pi_j\geq i$ for all $j>i$.  }
\smallskip

We consider the case $\pi_j\leq i$ for all $j<i$, the other case being similar. Let $L$ be a horizontal line drawn across $G_\pi$ at $y=i+\frac 12$. Since $G_\pi$ is connected, with components both above and below this horizontal line, it must contain vertical line segments crossing $L$.  These vertical line segments must correspond to indices $j$ where $j$ and $\pi_j$ lie on opposite sides of the value $i+\frac 12$.   Furthermore, as $\pi$ is a permutation, there must be an equal number of indices $j<i+\frac 12$ with $\pi_j>i+\frac 12$ as there are $j>i+\frac 12$ with $\pi_j<i+\frac 12$.

Given our hypothesis that $\pi_j\leq i$ for all $j<i$, there exists precisely one index $j$ that can satisfy the former property, namely $j=i$. Hence there exists also exactly one index $k>i$ with $\pi_k<i+\frac 12$.  Thus the diagram $G_\pi$ crosses $L$ twice, once at $i$ and again at $k$.  (Note that both of these line segments must be part of the outermost Seifert circle, $S_{\text{max}}$, and that $S_{\text{max}}$ therefore ``touches'' the diagonal at the point $(i,i)$.) We now produce a permutation $\pi'$, of length $n{+}1$, whose associated link corresponds to cutting the strands along the line $L$ and reattaching the two loose ends on either side. 

This new permutation $\pi'=\pi'_1\pi'_2\cdots\pi'_n\pi'_{n+1}$ is defined as follows. For $1\le j < i$ define $\pi'_j = \pi_j$; also, define $\pi'_i = \pi_k$ and $\pi'_{k+1} = i+1$. Finally, for $i\le j\le n$, $j\ne k$, set $\pi'_{j+1} = \pi_j+1$. Using that $G_\pi$ intersects $L$ only at the vertical lines in positions $i$ and $k$, it is not hard to check that $D(\pi') = D(\pi)$ and $I(\pi') = I(\pi)$. For example, the construction makes clear that if an inversion in $\pi$ is at elements that are both in positions less than $i$, or that both have positions greater than $i$, then there is a corresponding inversion in $\pi'$ (involving the same positions in the first case, and positions that are shifted by 1 in the latter case). An inversion in $\pi$ that involves positions $i$ and $j$, for $j>i$, corresponds to an inversion in $\pi'$ involving $i+1$ and $j+1$; in addition, inversions involving $j$, for some $j< i$, and position $k$ correspond to inversions in $\pi'$ at positions $j$ and $i$.

Since $|\pi'| =|\pi|+1$, we have that $x(\pi') = x(\pi) - 1$.

\begin{figure}[ht]
	\centering
    
     \begin{tikzpicture}[scale= 0.36,baseline={(0,0-0.15)}]
        \cycle{5,4,1,2,9,8,6,3,7}
        \draw[red,dash dot,line width=0.33mm] (-0.9,5.07)--(9.9,5.07);
        \draw [red] (10.05, 5.86) node {$L$};
      \end{tikzpicture}
	\hspace{0.8cm} \cyclefig{5,4,1,2,3,10,9,7,6,8}
	\caption{The cycle diagram of $\pi=541298637$ along with the line $L$ corresponding to the index $i=5$ satisfying the defining property of Case 2. 
 To the right is the permutation $\pi'$ obtained after ``cutting'' along $L$.}
	\label{fig:cutting}
\end{figure}
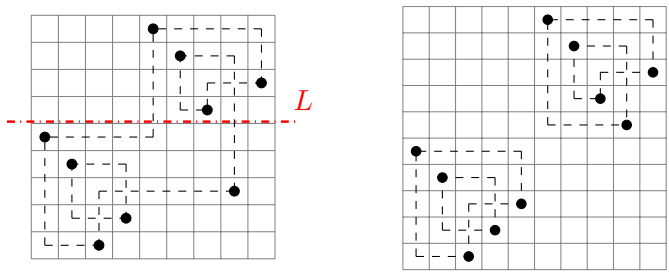

\begin{figure}[h!]
    \begin{tikzpicture}[>=stealth,scale=1.5]
    	\coordinate (A) at (0.3,1.6);
    	\coordinate (B) at (2.1,1.6);
    	\coordinate (C) at (0.3,1);
    	\coordinate (D) at (0.3,2.0);
    	\coordinate (E) at (2.1,0.9);
    	\coordinate (F) at (2.1,2.2);
    	\coordinate (K) at (0.723, 0.950);
    	\coordinate (L) at (0.484, 0.644);
    	\coordinate (M) at (1.133, 1.305);
    	\coordinate (N) at (1.365, 1.951);
    	\coordinate (O) at (0.824, 2.178);
    	\coordinate (P) at (0.639, 1.752);
    	\coordinate (Q) at (1.394, 1.386);
    	\coordinate (R) at (1.47, 0.890);
    	\coordinate (S) at (1.488, 2.222);
    	\coordinate (T) at (1.11, 0.416);

        \filldraw[gray!20!white] (0,0.1) -- (0,1) -- (0.3, 1) -- (0.3,2.6)--(2.1,2.6)--(2.1,0.1)--cycle;
        \draw[very thick] (0,0.1) -- (0,1) -- (0.3, 1) -- (0.3,2.6);
       	\draw[very thick] (2.1,0.1) -- (2.1,2.6);
       	
       	\draw (D) -- (O) -- (P) -- (D);
       	\draw (O) -- (M) -- (P);
       	\draw (O) -- (S) -- (N) -- (F) -- (S);
       	\draw (O) -- (N) -- (M) -- (Q) -- (N) -- (B) -- (Q);
       	\draw (M) -- (K) -- (P) -- (A) -- (K) -- (C) -- (L) -- (K) -- (T) -- (L);
       	\draw (T) -- (M) -- (R) -- (T) -- (E) -- (R) -- (Q);
       	\draw (B) -- (R);
       	
       	\draw 
       		($(L) + (-0.3,-0.4)$) -- (L) -- ($(L) + (0.15,-0.4)$)
       		($(T) + (-0.2,-0.3)$) -- (T) -- ($(T) + (0.5,-0.25)$)
       		($(S) + (-0.06,0.3)$) -- (S) -- ($(S) + (0.15,0.3)$)
       		($(O) + (-0.18,0.3)$) -- (O) -- ($(O) + (0.25,0.3)$);
       	\foreach \p in {(K),(L),(M),(N),(O),(P),(Q),(R),(S),(T)}
			\filldraw \p circle (1pt);
       	\filldraw[red] 
       		(A) circle (2pt)
       		(B) circle (2pt);
       	\filldraw
       		(C) circle (2pt)
       		(D) circle (2pt)
       		(E) circle (2pt)
       		(F) circle (2pt);
       	\draw[red] (A) -- (B);
       	\draw[red] 
       		(A) node[left]{$A$}
       		(B) node[right]{$B$};
       	
       	\filldraw[gray!20!white] (3,0.1) -- (3,1) -- (3.3, 1) -- (3.3,2.6)--(5.1,2.6)--(5.1,0.1)--cycle;
       	\draw[very thick] (0+3,0.1) -- (0+3,1) -- (0.3+3, 1) -- (0.3+3,2.6);
       	\draw[very thick] (2.1+3,0.1) -- (2.1+3,2.6);
       	
       	\draw ($(D)+(3,0)$) -- ($(O)+(3,0)$) -- ($(P)+(3,0)$) -- ($(D)+(3,0)$);
       	\draw ($(O)+(3,0)$) -- ($(S)+(3,0)$) -- ($(N)+(3,0)$) -- ($(F)+(3,0)$) -- ($(S)+(3,0)$);
       	\draw ($(O)+(3,0)$) -- ($(N)+(3,0)$);
       	\draw ($(M)+(3,0)$) -- ($(Q)+(3,0)$);
       	\draw ($(N)+(3,0)$) -- ($(B)+(3,0)$) -- ($(Q)+(3,0)$);
       	\draw ($(M)+(3,0)$) -- ($(K)+(3,0)$) -- ($(A)+(3,0)$) -- ($(P)+(3,0)$);
       	\draw ($(K)+(3,0)$) -- ($(C)+(3,0)$) -- ($(L)+(3,0)$) -- ($(K)+(3,0)$) -- ($(T)+(3,0)$) -- ($(L)+(3,0)$);
       	\draw ($(T)+(3,0)$) -- ($(M)+(3,0)$) -- ($(R)+(3,0)$) -- ($(T)+(3,0)$) -- ($(E)+(3,0)$) -- ($(R)+(3,0)$) -- ($(Q)+(3,0)$);
       	\draw ($(B)+(3,0)$) -- ($(R)+(3,0)$);
       	\draw 
       	($(L) + (3-0.3,-0.4)$) -- ($(L)+(3,0)$) -- ($(L) + (3.15,-0.4)$)
       	($(T) + (3-0.2,-0.3)$) -- ($(T)+(3,0)$) -- ($(T) + (3.5,-0.25)$)
       	($(S) + (3-0.06,0.3)$) -- ($(S)+(3,0)$) -- ($(S) + (3.15,0.3)$)
       	($(O) + (3-0.18,0.3)$) -- ($(O)+(3,0)$) -- ($(O) + (3.25,0.3)$);
       	
       	\coordinate (AA) at ($(A) + (3,0)$);
       	\coordinate (BB) at ($(B)+(3,0)$);
       	\coordinate (KK) at ($(K)+(3,0)$);
       	\coordinate (LL) at ($(L)+(3,0)$);
       	\coordinate (MM) at ($(M)+(3,0)$);
       	\coordinate (NN) at ($(N)+(3,0)$);
       	\coordinate (OO) at ($(O)+(3,0)$);
       	\coordinate (PP) at ($(P)+(3,0)$);
       	\coordinate (QQ) at ($(Q)+(3,0)$);
       	\coordinate (RR) at ($(R)+(3,0)$);
       	\coordinate (SS) at ($(S)+(3,0)$);
       	\coordinate (TT) at ($(T)+(3,0)$);
       	
       	\draw[red] 
       		(QQ) -- (AA) -- (MM);
       	\draw[red]
       		(AA) ..controls ($(AA)!0.2!(BB) + (0,0.08)$) and ($(AA)!0.4!(BB)$) ..(OO)
       		(BB) ..controls ($(AA)!0.8!(BB) + (0,0.08)$) and ($(AA)!0.6!(BB)$) ..(OO);
       	\foreach \p in {(KK),(LL),(MM),(NN),(OO),(PP),(QQ),(RR),(SS),(TT)}
       	\filldraw \p circle (1pt);
       	\filldraw[red] 
       	($(A)+(3,0)$) circle (2pt)
       	($(B)+(3,0)$) circle (2pt);
       	\filldraw
       	($(C)+(3,0)$) circle (2pt)
       	($(D)+(3,0)$) circle (2pt)
       	($(E)+(3,0)$) circle (2pt)
       	($(F)+(3,0)$) circle (2pt);
       	\draw[red] ($(A)+(3,0)$) -- ($(B)+(3,0)$);
       	\draw[red] 
       	($(A)+(3,0)$) node[left]{$A$}
       	($(B)+(3,0)$) node[right]{$B$};
    \end{tikzpicture}
    \caption{Cutting $F_{\pi}$ along $L$: Let $A, B$ be the points where $L$ intersects the boundary of the surface. A subdivision of a given triangulation will have $A$ and $B$ as vertices. Remove all edges that intersect $L$, and vertices (other than $A$ and $B$) that intersect $L$ as well as their incident edges. What remains is a continuous image of a polygon with $L$ as a secant from one vertex to another. Triangulate this polygon with edges emanating only from $A$ and $B$.}
    \label{fig:cut_along_L}
\end{figure}
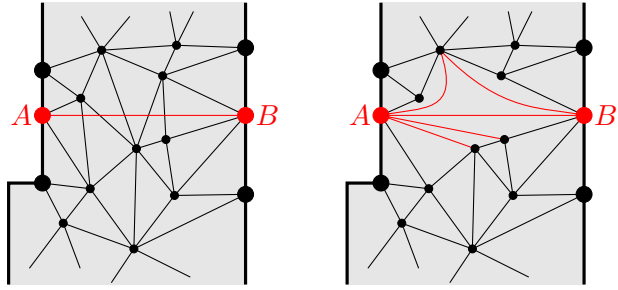

Recall that part of the surface $F_\pi$ is a disk in a horizontal plane of $\mathbb R^3$, having $S_{\text{max}}$ as its boundary. Viewing $L$ as being contained in that plane, the surface $F_{\pi'}$ is homeomorphic to the result of cutting $F_\pi$ along the line $L$. Consequently, the Euler characteristic changes as $\chi(F_{\pi'}) = \chi(F_\pi) + 1$ (there is a triangulation so that the cut is along one edge between two vertices on the boundary of the surface, cf.\ Figure \ref{fig:cut_along_L}; cutting can be seen to duplicate the edge and the two vertices). Thus, the desired equality holds for $\pi$ if and only if it holds for $\pi'$. 

Since $G_{\pi'}$ is not connected, the argument from Case 1 allows us to decompose $\pi'$ as $\rho\oplus\tau$. Additionally, $\pi'$ cannot have a fixed point. This implies that the length of $\rho$ and the length of $\tau$ are at most two less than $|\pi'|$; that is, $\max\{|\rho|, |\tau|\} \le n-1$. 

\subsection{Case 3.} \textit{The diagram $G_\pi$ is connected and no index with the property described in Case 2 exists.}
\smallskip

In this case, we construct a new permutation $\pi'$ with $|\pi'| = |\pi|-2$.  As we will see, the key observation is that all of the line segments occurring in the cycle diagram of this new permutation $\pi'$ are also present in the cycle diagram of $\pi$, but all of those line segments in $G_\pi$ corresponding to the outermost Seifert circle have been removed.  Furthermore, we will see that $\pi'$ has the property that $x(\pi')=-\chi(F_{\pi'})$ if and only if $x(\pi)=-\chi(F_{\pi})$.  

Some preliminaries are needed in order to define $\pi'$ and prove that it has the desired property. First, as in Case 2, the connected diagram $G_\pi$ has a unique maximal Seifert circle $S_{max}$. To define $\pi'$ we need first to identify the columns in the diagram where $S_{max}$ has vertical line segments. These columns will necessarily be indexed by $1$, $n$, and every column in which $G_\pi$ has a crossing that creates a part of the Seifert circle $S_{max}$. 
Say that there are $a$ such crossings that are above the diagonal and $b$ such crossings that are below the diagonal.

Write $i_1 < i_2 < \cdots < i_{a}$ for those columns where $S_{max}$ has a crossing above the diagonal. On the other side, write $j_1 > j_2 > \cdots > j_{b}$ be the columns where $S_{max}$ has a crossing below the diagonal, written in decreasing order. Finally, define $i_0 = 1 = j_{b+1}$ and $i_{a+1}=n=j_0$.  

Notice, by the hypothesis of Case 3 that no index with the property described in Case 2 exists, it is necessary that the outermost Seifert circle $S_{max}$ touches the main diagonal only at $(1,1)$ and $(n,n)$, and so every vertical line segment of $S_{max}$ is in one of the columns indexed above.    Thus, it must be that $\pi_{i_{a}} = n$ and $\pi_{j_{b}} = 1$. 

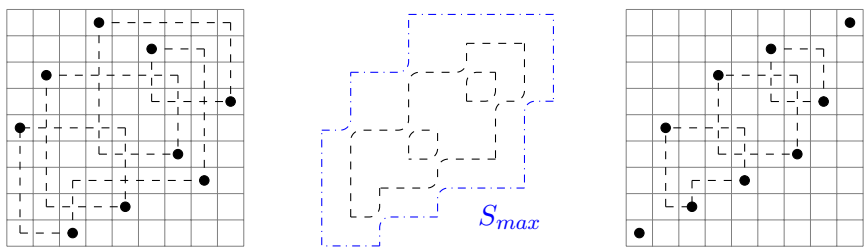
\begin{figure}[ht]
	\centering
	\ \cyclefig{5,7,1,9,2,8,4,3,6}
	\hspace{0.5cm} \begin{tikzpicture}[scale=3.85]
	\draw[blue,dash dot] (0.3,0.1)--(0.1,0.1)--(0.1,0.5);
	\draw [rounded corners, blue,dash dot] (0.1,0.5)--(0.2,0.5)--(0.2,0.7);
	\draw [rounded corners, blue,dash dot] (0.2,0.7)--(0.4,0.7)--(0.4,0.9);
	\draw[blue,dash dot] (0.4,0.9)--(0.9,0.9)--(0.9,0.6);
	\draw [rounded corners, blue,dash dot] (0.9,0.6)--(0.8,0.6)--(0.8,0.3);
	\draw[rounded corners, blue,dash dot] (0.8,0.3)--(0.5,0.3)--(0.5,0.2);
	\draw[rounded corners, blue,dash dot] (0.5,0.2)--(0.3,0.2)--(0.3,0.1);
	\draw[rounded corners,dashed] (0.2,0.2)--(0.2,0.5)--(0.4,0.5)--(0.4,0.7)--(0.6,0.7)--(0.6,0.8);
	\draw [dashed] (0.6,0.8)--(0.8,0.8);
	\draw[rounded corners,dashed] (0.8,0.8)--(0.8,0.6)--(0.7,0.6)--(0.7,0.4);
	\draw [rounded corners,dashed] (0.7,0.4)--(0.5,0.4)--(0.5,0.3)--(0.3,0.3);
	\draw [rounded corners,dashed] (0.3,0.3)--(0.3,0.2)--(0.2,0.2);
	\draw[rounded corners,dashed] (0.6,0.6)--(0.6,0.7)--(0.7,0.7);
	\draw [rounded corners,dashed] (0.7,0.7)--(0.7,0.6)--(0.6,0.6);
	\draw[rounded corners,dashed] (0.4,0.42)--(0.4,0.5)--(0.5,0.5);
	\draw [rounded corners,dashed] (0.5,0.49)--(0.5,0.4)--(0.4,0.4);
	\draw [blue] (0.75, 0.2) node {$S_{max}$};
	\end{tikzpicture}
	\hspace{0.7cm}
	\cyclefig{1,5,2,7,3,8,4,6,9}
	
	\caption{The cycle diagram of $\pi =  571928436$ and its corresponding Seifert circles. The Seifert circle $S_{max}$ is labeled and drawn in blue (and dash-dotted). In this case we have $a=2$, $b=3$, with $i_0=1$,  $i_1=2$, $i_2=4$ and $j_0=9$, $j_1=8$, $j_2=5$, $j_3=3$. To the right is the permutation $\bar{\pi}=152738469 = \pi \circ (4\ 2\ 1\ 3\ 5\ 8\ 9)$.}
	\label{fig:max_Seifert}
\end{figure}

Now, we define an intermediate permutation $\bar \pi = \pi \circ (i_a \ i_{a-1} \cdots  i_1\ 1\ j_b \ j_{b-1} \cdots j_1\ n)$, using cycle notation.  Alternatively, the construction of $\bar\pi$ can be made as follows. Beginning from the one-line notation $\pi=\pi_1\pi_2\ldots\pi_n$: move $\pi_1$ into position $i_1$, then move $\pi_{i_1}$ into position $i_2$. Continue this way until moving $\pi_{i_a}=n$ into position $i_{a+1}=n$. Next, move $\pi_n = \pi_{j_0}$ into position $j_1$, move $\pi_{j_1}$ into position $j_2$, and so on. Finally, move $\pi_{j_b}=1$ into position $j_{b+1}=1$. No change is made to other elements in the one-line notation; the resulting permutation is $\bar\pi$.

The one-line notation $\bar\pi = \bar\pi_1\bar\pi_2\ldots\bar\pi_n$ is the following. 
\begin{equation}
    \bar\pi_{\ell} = \begin{cases}
    \pi_{i_{k-1}} & \text{if $\ell = i_k$ for some $1\le k\le a+1$}\\
     \pi_{j_{k-1}} & \text{if $\ell = j_k$ for some $1\le k\le b+1$}\\
      \pi_\ell & \text{otherwise}.
    \end{cases}
\end{equation}
Note that $|\bar\pi| = |\pi|$ and that $\bar\pi$ fixes both $1$ and $n$. The permutation $\pi'$, used in the induction, is the permutation such that $\bar\pi = \text{id}_1 \oplus \pi' \oplus \text{id}_1$ (where $\text{id}_1$ is the identity permutation on one element), and can be obtained from $\bar\pi$ by setting $\pi'_k = \bar\pi_{k+1}-1$ for $1\le k\le n-2$.  
Note, the link type of $L_{\pi'}$ is generally different than the link type of $L_{\pi}$. 

As mentioned above,  
we have the following relationship between the Seifert circles of $G_\pi$ and those of $G_{\pi'}$, the verification of which is left to the reader.

\begin{lemma}The diagram $G_{\pi'}$ is the result of removing $S_{max}$ from $G_\pi$ (after only smoothing the crossings on $S_{max}$).
\end{lemma}

Setting $c = a+b$, so that $c$ is the number of crossings of $G_{\pi}$ on $S_{max}$, we see that $G_{\pi'}$ has exactly $c$ fewer crossings than $G_{\pi}$. It also has exactly one less upper right corner, the one from $S_{max}$. Since $-\chi(F_{\pi})$ may be determined \eqref{eq:chi} as the number of crossings of $G_{\pi}$ minus the number of upper right corners of $G_{\pi}$ (and likewise for $\pi'$), we get that \[-\chi(F_{\pi}) + \chi(F_{\pi'}) = c - 1.\]

We may finish Case 3, and so complete the proof of Theorem \ref{thm:main}, by showing that $x(\pi) - x(\pi') = c-1$. The remainder of the paper is devoted to doing just this, with the bulk of that effort focused on understanding how inversions in $\pi$ compare to those in $\pi'$.

To begin, we consider certain sets of pairs of elements from $\pi$. Given $0\le k \le a$, define $A_k = \{(\pi_\ell, \pi_{i_k})\ :\ i_k<\ell<i_{k+1}\}$. Given $0\le k \le b$, define $B_k = \{(\pi_{j_k}, \bar\pi_\ell)\ :\ j_{k+1}<\ell<j_{k}\}$, with the one exception that the pair $(\pi_n, n)$ is also in $B_{0}$ (i.e., noting that $j_{0} = n = \bar\pi_n$, make the upper bound on $\ell$ be a weak upper bound for $k=0$). Finally, define \[A = \bigcup_{k=0}^{a}A_k;\qquad B =\bigcup_{k=0}^{b}B_k.\]

Regarding the definition of $B_k$ note that if $j_{k+1} < \ell < j_k$, with $\ell \ne i_{k'}$ for any $k'$, then $\bar\pi_\ell = \pi_\ell$. If it is the case that $\ell = i_{k'}$ for some $1\le k'\le a$ then $\bar\pi_\ell = \pi_{i_{k'-1}}$.

\begin{lemma}\label{lem:inversions}
Every ordered pair $(p,q) \in  A\cup B$ is an inversion in $\pi$.
\end{lemma}
\begin{proof} Let $k$ and $\ell$ be such that $0\le k\le a$ and $i_{k} < \ell < i_{k+1}$. In the case that $k=a$, $\pi_{i_k} = n$, and so it is clear that $(\pi_\ell, \pi_{i_k})$ is an inversion in $\pi$.

Now, for $k$ with $0\le k< a$, there is a crossing of $G_\pi$ in the crossings corresponding to $S_{max}$ at grid coordinates $(i_{k+1}, \pi_{i_{k}})$. If it were the case that 
$\pi_{i_k} < \pi_\ell$ for some $\ell$ with $i_k < \ell < i_{k+1}$ then the vertical line segment extending from index $\ell$ would create a crossing at coordinates $(\ell,\pi_{i_k})$. But then a Seifert circle corresponding to this crossing would contain a part of $S_{max}$ in its interior, which contradicts $S_{max}$ being maximal. Therefore, every pair in $A$ is an inversion in $\pi$.

The pair $(\pi_n, n) \in B_0$ is an inversion in $\pi$, as $G_\pi$ being connected requires that $\pi_n < n$, and $(\pi_n, n) = (\pi_n, \pi_{i_a})$.  Now, suppose that $\ell$ is such that, for some $0\le k\le b$ we have $j_{k+1} < \ell < j_k$ and that $\ell \ne i_{k'}$ for all $1\le k'\le a$.  Since $\bar\pi_\ell = \pi_\ell$ in this case, a pair $(\pi_{j_k},\bar\pi_\ell)\in B$ is equal to $(\pi_{j_k}, \pi_\ell)$. That this is an inversion in $\pi$ follows from an analogous argument to the previous one for pairs in $A$. 

Finally, suppose for some $0\le k\le b$ there exists such a $k'$, so that $j_{k+1} < i_{k'} < j_k$. Recall that $\pi_{i_{k'-1}}=\bar\pi_{i_{k'}}$. We must have that $\bar\pi_{i_{k'}} > i_{k'}$ since there is a crossing above the diagonal at coordinates $(i_{k'}, \pi_{i_{k'-1}})$. Additionally, considering the crossing corresponding to column $j_{k+1}$, which is below the diagonal (or, that $\pi_{j_b}=1=j_{b+1}$), we have $\pi_{j_k} \le j_{k+1}$. Hence, $\pi_{j_k} \le j_{k+1} < i_{k'} < \pi_{i_{k'-1}}$.  Since $i_{k'-1}< i_{k'} < j_k$, we see that $(\pi_{j_k}, \bar\pi_{i_{k'}}) = (\pi_{j_k}, \pi_{i_{k'-1}})$ is an inversion in $\pi$.
\end{proof}

The definition of the sets $B_k$, particularly the use of pairs having the form $(\pi_{j_k}, \bar\pi_\ell)$ for $j_{k+1}<\ell<j_k$, rather than the form $(\pi_{j_k}, \pi_\ell)$, is convenient. On the one hand, it ensures that we find certain inversions of $\pi$ in $A\cup B$ {--} namely, pairs $(\pi_{j_k}, \pi_{i_{k'-1}})$ with $j_{k+1} < i_{k'} < j_k$ that may exist. In addition, the definition helps by making $A$ and $B$ disjoint.

\begin{lemma}\label{lem:disjoint}
    Let $0\le k\le a$ and $0\le k'\le b$. Then $A_k \cap B_{k'} = \emptyset$.
\end{lemma}
\begin{proof}
    Fix $0\le k\le a$ and $0\le k'\le b$ and suppose $(p,q) \in A_k \cap B_{k'}$. We have $p = \pi_{j_{k'}}$ and $q = \pi_{i_k}$. By the definition of $A_k$, we have $i_k < j_{k'} < i_{k+1}$ (preventing $k'$ from being 0). In addition, $q = \bar\pi_\ell$ for some $j_{k'+1} < \ell < j_{k'}$. However, $q = \pi_{i_k}= \bar\pi_{i_{k+1}}$, and so $\ell =  i_{k+1}$, a contradiction. And so $A_k \cap B_{k'} = \emptyset$.
\end{proof}
We note that, from their definition, it is clear that the sets $A_k$, $0\le k\le a$ are themselves pairwise disjoint; as are the sets $B_k$, $k=0,1,\ldots,b$.

\begin{lemma}\label{lem:inversion-locations}
  Fix $p,q$ with $1\leq p < q \leq |\pi|$. Then the relative position of $p$ and $q$ in $\pi$ is different than the relative position of $p$ and $q$ in $\bar\pi$ if and only if
	$$(p,q) \in A\cup B.$$
\end{lemma}
Note that just as $p < q$ in the statement, the definitions of $A_k$ and $B_k$ are such that ordered pairs in these sets have a first coordinate that is less than the second coordinate.

\begin{proof}
We have already observed that each $(p,q)\in A\cup B$ is an inversion in $\pi$, meaning that $q$ appears to the left of $p$ in $\pi$. If $(p,q) \in B_k$ for some $k$, with $q=\bar\pi_\ell$ for some $j_{k+1}<\ell<j_k$, then $q$ appears to the right of $p=\pi_{j_k}$ in $\bar\pi$ since $\bar\pi_{j_{k+1}}=\pi_{j_k}$. Now, if $(p,q)\in A_k$ for some $0\le k\le a$, then $\pi_{i_k}$ appears in position $i_{k+1}$ in $\bar\pi$. Letting $\ell$ be such that $i_k < \ell < i_{k+1}$, then $\pi_\ell$ will either appear in position $\ell$ in $\bar\pi$ or it will appear in a position that is farther left (when $\ell = j_{k'}$, since $j_{k'+1} < j_{k'}$). 

In both cases, $q = \pi_{i_k}$ appears to the right of $p=\pi_\ell$ in $\bar\pi$.  Thus, we may conclude that the relative position of $p$ and $q$ does change if $(p,q) \in A\cup B$.

Suppose that $p < q$ are elements such that the relative position of $p$ and $q$ in $\pi$ is different than their relative position in $\bar\pi$. It must be that one of $p$ and $q$, at least, has position in $\pi$ that is in $\{i_0,i_1,\ldots,i_a\} \cup \{j_0,j_1,\ldots,j_b\}$  (since elements in other positions in $\pi_1\pi_2\ldots\pi_n$ do not change their position). 

In the first case, let $x = \pi_{i_k}$ for some $0\le k\le a$, and that $x$ equals one of $p$, $q$. Let $y$ be the other element of $\{p,q\}$. Considering the construction of $\bar\pi$, $x$ appears in $\bar\pi$ farther to the right than its position in $\pi$. For the relative position of $x$ and $y$ to change, $y$ cannot be in position $i_{k'}$ for some $0\le k'\le a$, since $\pi_{i_k}$ and $\pi_{i_{k'}}$ have the same relative position in $\pi$ and $\bar\pi$. Hence, $y$ must appear in $\bar\pi$ farther to the left than it appears in $\pi$, or in the same position as in $\pi$.
As a consequence, in $\pi$ we have $y$ appearing to the right of $x$; moreover, its position is either left of $i_{k+1}$ or equal to $j_{k'}$ where $j_{k'+1} < i_{k+1} \le j_{k'}$ (the weak inequality needed when $k'=0$). In the first case, Lemma \ref{lem:inversions} and the fact that $p < q$ together imply that $q = x = \pi_{i_k}$ and $p = \pi_\ell$ for some $i_k<\ell<i_{k+1}$; we find that $(p,q) \in A$. In the second case we have that $j_{k'+1} < i_{k+1} \le j_{k'}$ and $x = \pi_{i_k} = \bar\pi_{i_{k+1}}$. By Lemma \ref{lem:inversions} (which implies that $\pi_{j_{k'}} < \bar\pi_{i_{k+1}}$), it must be that $x = q$ and we conclude that $(p,q) \in B$.

The remaining case is that one of $p,q$ is equal to $\pi_{j_k}$ for some $0\le k\le b$, and that the other is $\pi_\ell$ for some $\ell$, and is such that $\pi_{\ell}=\bar\pi_{\ell}$.  Since the relative position of $p, q$ is different in $\bar\pi$ compared to in $\pi$, and $\pi_{j_k}=\bar\pi_{j_{k+1}}$, it must be that $j_{k+1} < \ell < j_k$. As shown in Lemma \ref{lem:inversions}, we must have $\pi_{j_k} < \pi_\ell$, so $p = \pi_{j_k}$ and $q = \pi_\ell = \bar\pi_\ell$. And so $(p, q) \in B$, concluding the proof of the lemma.
\end{proof}

Let $(p,q)$ be an ordered pair with $p<q$. By combining Lemma \ref{lem:inversions} and Lemma \ref{lem:inversion-locations}, if $(p,q)$ is not an inversion in $\pi$ then it cannot be an inversion in $\bar\pi$. We see that the cardinality of $A\cup B$ equals $I(\pi) - I(\bar\pi)$.

Let $c = a+b$ be the number of crossings in $S_{max}$. Using Lemma \ref{lem:disjoint} we have
\begin{align*}
    I(\pi) - I(\bar\pi)   &= \sum_{k=0}^{a}|A_k| + \sum_{k=0}^{b}|B_k| \\
    	&= \sum_{k=0}^{a}(i_{k+1} - i_k - 1) + (j_0 - j_1) + \sum_{k=1}^{b}(j_k - j_{k+1} - 1) \\
            &= (i_{a+1} - i_0) + (j_0 - j_{b+1}) - (c+1)  \\
            &= 2n-3 - c,
\end{align*}
the last equation since $i_{a+1} = n = j_0$, $j_{b+1} = 1 = i_0$.

The total displacement of the permutation changes as follows. For each $i_k$, $1\le k\le a$, recall $\bar\pi_{i_k} = \pi_{i_{k-1}}$. Since $\pi_{i_{k-1}} < \pi_{i_k}$ and also $i_k$ is less than both $\pi_{i_k}$ and $\pi_{i_{k-1}}$, we have $|\pi_{i_k} - i_k| - |\bar\pi_{i_k} - i_k| = \pi_{i_k} - \pi_{i_{k-1}}$. In addition, $|\pi_{i_0} - i_0| - |\bar\pi_{i_0} - i_0| = \pi_{i_0} - 1$, since $\bar\pi_{i_0} = i_0=1$. 

Similarly, we have $|j_k - \pi_{j_k}| - |j_k - \bar\pi_{j_k}| = \pi_{j_{k-1}} - \pi_{j_{k}}$ for each $1\le k\le b$, and we have $|\pi_{j_0} - j_0| - |\bar\pi_{j_0} - j_0| = n - \pi_n$. Putting this together,
\[D(\pi) - D(\bar\pi) = \pi_{i_{a}} - 1 + n - \pi_{j_{b}} = 2(n - 1).\]
Recall the permutation $\pi'$, defined so that $\bar\pi = \text{id}_1 \oplus \pi' \oplus \text{id}_1$. It is clear that $D(\bar\pi) = D(\pi')$, $I(\bar\pi) = I(\pi')$, and $|\bar\pi| - |\pi'| = 2$. Therefore,
\begin{align*}
    x(\pi) - x(\pi') = 2(n-1) - (2n-3 - c) - 2 = c - 1,
\end{align*}
proving what was needed for the induction step in Case 3.

\bibliographystyle{amsplain}
\bibliography{ldgbib}

\end{document}